\documentclass{ifacconf}

\usepackage{graphicx}      
\graphicspath{{./}{Figure/}}
\usepackage{natbib}        

\usepackage{amsmath}
\usepackage{amssymb}
\usepackage{mathrsfs,graphicx,paralist}
\usepackage{color}
\theoremstyle{plain}\newtheorem{remark}{Remark}[section]
\theoremstyle{plain}

\newcommand{\R}{\mathbb{R}}

\def\T{\mathcal{T}}

\begin{document}
\begin{frontmatter}

\title{The Hughes model for pedestrian dynamics and congestion modelling\thanksref{footnoteinfo}} 

\thanks[footnoteinfo]{EC  acknowledges financial support from INDAM GnCS, project ``Metodi numerici semi-impliciti e semi-Lagrangiani per sistemi iperbolici di leggi di bilancio''.
AF acknowledge financial support from the Austrian Academy of Sciences \"OAW via the New Frontiers Group NST-001.     FJS benefited from the support of the ``FMJH Program Gaspard Monge in optimization and operation research'', and from the support to this program from EDF.}

\author[Carlini]{Elisabetta Carlini} 
\author[Festa]{Adriano Festa} 
\author[Silva]{Francisco J. Silva}

\address[Carlini]{ Dipartimento di Matematica ``G. Castelnuovo'', Sapienza Universit\`a di Roma, (e-mail: carlini@mat.uniroma1.it).}
\address[Festa]{ RICAM -- Johann Radon Institute
for Computational and Applied Mathematics (\"OAW),  (e-mail: adriano.festa@oaew.ac.at)}
\address[Silva]{ XLIM - DMI 
UMR CNRS 7252 Facult\'e des Sciences et Techniques, Universit\'e de Limoges,   (e-mail: francisco.silva@unilim.fr)}

\begin{abstract}                
In this paper we present  a numerical study  of some variations of the Hughes model for pedestrian flow under different types  of  congestion effects. The  general model consists of a coupled non-linear PDE system involving an eikonal equation and a first order conservation law, and it intends to approximate the flow of a large pedestrian group aiming to reach a target as fast as possible, while taking into account the congestion of the crowd.   

We propose an efficient semi-Lagrangian scheme (SL) to approximate the solution of the PDE system    and  we investigate the macroscopic effects of different penalization functions  modelling the congestion  phenomena. 
\end{abstract}


\end{frontmatter}

\section{Introduction}

In the last years the understanding of the dynamics underneath the   pedestrian crowd motion has attracted a wide attention in public debate and in the scientific community. The complete comprehension of such phenomena remains still as a subject of active research. Nevertheless, important progresses  have been made in the past two decades and several mathematical models have been proposed for its description. \\ 
\noindent The different approaches to tackle the problem can be classified according to the scale of the model: they range from microscopic systems, where each  trajectory is described individually, to macroscopic systems, where   the dynamics of the crowd is  modelled as a time-evolving  density distribution of indistinguishable agents.
 For a detailed overview we refer e.g. to the work  by \cite{bellomo:2011} and the monograph by \cite{cristiani2014multiscale}.
\\
In the article by \cite{h00}, the author introduced a by now classical  fluid-dynamical  model   to study the motion  of a large human crowd. The crowd is treated as a ``thinking fluid'' and it moves at  maximum speed towards  a common target, while  taking  into account   environmental factors such as the congestion of the crowd. In  this model people prefer  to avoid crowded regions and this assumption is incorporated in a  potential  function  which determines the velocity field driving the crowd. Indeed, the potential is characterized as  the  gradient of the solution of an eikonal equation with a right hand side which depends on the  current distribution of the crowd.
Hence,   given a time horizon $T>0$, for each time instant $t\in [0,T]$, each microscopical individual  looks at the global configuration of the crowd and updates his/her velocity trying to avoid the  current  crowded regions. \\
\noindent In a two-dimensional space setting, given $f: [0,+\infty[ \to \R$, a bounded   domain $\Omega\in \R^2$ and a target $\mathcal{T}\subseteq \Omega$, the model is given by
\begin{equation}\label{hughes}
\left\{
\begin{array}{ll}
\partial_t m(x,t)-\hbox{div}(f^2(m(x,t))\nabla u(x,t)\,m(x,t))=0,  \\[5pt]
\hspace{3.9cm}\mbox{in $ \Omega \times ]0,T[$,}\\[4pt]
|\nabla u(x,t)|=\displaystyle\frac{1}{f(m(x,t))},  \; \; \mbox{in } \left(\Omega \setminus \mathcal{T}\right) \times ]0,T[,
\end{array}\right.
\end{equation}
 complemented with the boundary conditions
 \begin{equation}\label{BC}
\hspace{-0.18cm}\left\{
\begin{array}{ll}
m(x,0)= m_0(t),  &\text{ on  } \Omega\times\{0\},\\[4pt]
u(x,t) = 0,  & \text{ on } \T \times (0,T),\\[4pt]
u(x,t)= g(x)   &\text{ on }   { \partial  \Omega } \times (0,T),  \\[4pt]
(f^2(m)\nabla u \,m)(x,t)\cdot \hat{n} (x)= 0,  &\text{ on } \partial \Omega \times (0,T).
\end{array}\right.
\end{equation}
In \eqref{hughes} the gradient operator $\nabla$ acts on the space variable $x$, the unknown variable $m$ is the pedestrian density and $u$  is the potential, i.e. the weighted shortest distance to the target $\mathcal{T}$. Concerning the boundary conditions \eqref{BC}, we assume that  $g$ is a continuous function large enough to satisfy some compatibility conditions,  see e.g. \cite[Chapter IV]{BardiCapuz@incollection}, and $\hat n$ denotes the outward normal to the boundary $\partial \Omega$. 
System \eqref{hughes} is a highly non-linear coupled system of PDEs. Few analytic results are available, all of them restricted to spatial dimension one and particular choices for the function $f$ (see the works by  \cite{di2011hughes} and \cite{amadori2014existence}). The main difficulty comes from the low regularity of the potential $u$, which is only Lipschitz-continuous.

Hughes proposed a few functions $f$'s penalizing regions of high density, the simplest choice being $f(m)=1-m$ where $1$ corresponds to the maximum scaled pedestrian density. In this work we focus our attention on numerical methods to solve \eqref{hughes}-\eqref{BC} for several choices of the penalty function $f$. We must underline that since well-posedness results for \eqref{hughes}-\eqref{BC} have not been proved  in full generality yet, convergence results of numerical algorithms for   \eqref{hughes}-\eqref{BC} seem currently out of reach. Thus, we will consider some heuristics to solve \eqref{hughes}-\eqref{BC} based on recent techniques introduced in \cite{CS12}, \cite{CS15} and \cite{carlini2016PREPRINT}. We point out that the article by \cite{carlini2016PREPRINT} concerns the approximation of a regularized version of \eqref{hughes}-\eqref{BC}. Our strategy is to look at \eqref{hughes}-\eqref{BC} as a single non-linear continuity equation, which can be approximated by discretizing    the characteristics governing the equation with  an Euler scheme. However,  the velocity field describing the characteristics depends at each time step non-locally on the current distribution of the agents and has to be computed with the help of the eikonal equation, for which several efficient methods exist.



\section{A semi-Lagrangian scheme for the approximation of the system}\label{discretization}
\subsection{Modeling aspects}\label{mod} As in the theory of Mean Field Games (MFGs), system \eqref{hughes} can, {\it at least formally}, be interpreted as a sort of Nash equilibrium for a dynamical game involving a continuous number of agents. Indeed, given a space-time distribution density  of the agents $(x,t) \in \Omega \times  [0,T] \mapsto m(x,t) \in \R$, such that $m \geq 0$, the solution $u[m]$ of the second equation in  \eqref{hughes}   can be formally represented as   the value function of the \emph{optimal control problem} 
\begin{equation}\label{value_function_depends_on_m}
u[m](x,t)= \inf_{ \alpha \in \mathcal{A}} \int_{t}^{\tau^{x,t}[\alpha]} F(m(X^{x,t}[\alpha](s),t)) d s, 
\end{equation}
where 
$$\begin{array}{l}
\mathcal{A} :=  \left\{ \alpha :[0,T] \mapsto \R^{d} \; ; \; |\alpha(t)| \leq 1, \; \mbox{a.e. in $[0,T]$} \right\},\\[6pt]
X^{x,t}[\alpha](s) :=  x+ \int_{t}^{s} \alpha(r) dr \hspace{0.4cm} \forall \; s\in [t,T],\\[7pt]
\tau^{x,t}[\alpha] :=  \inf\left\{ s\in [t,T] \; ; \; X^{x,t}[\alpha](s) \in \mathcal{T}\right\},\\[6pt]
 \mbox{and }  \;  F(m):= \frac{1}{f(m)}.
\end{array}$$
Thus, $u[m](x,t)$  corresponds to the {\it weighted minimal time}  to reach the target set $\mathcal{T}$
 for a   {\it typical player} positioned at $x$ at time $t$. From the modelling point of view, it is fundamental to notice that in its individual cost  the typical agent {\it freezes} the global distribution $m$ at time $t$ and thus   s/he does not forecast  the future behaviour of the population in order to determine its optimal policy. This modelling aspect shows an important difference with MFG models, where, at the equilibrium,  the agents take into account the future distribution of the population in order to design their strategies. 
 
Then, after the optimal feedback law
\begin{equation}\label{feedbacklaw} s\in [t,T] \mapsto  \bar{\alpha}[t] (s):= -\frac{\nabla u[m](X^{x,t}[\hat{\alpha}](s),t)}{\left|\nabla u[m](X^{x,t}[\hat{\alpha}](s),t)\right|}\end{equation}
is computed, the agents {\it actually} move according to the dynamics  defined by the the solution of the ODE 
\begin{equation}\label{evolucionconfeedback}
\frac{d \hat{X}^{t,x}}{d s}(s)=  -\frac{\nabla u[m](\hat{X}^{x,t}(s),s)}{\left|\nabla u[m](\hat{X}^{x,t}(s),s)\right|} f(m(\hat{X}^{t,x}(s),s)),
\end{equation}
for $s\in [t,T]$. 
Note that at each time $s \in [s,T]$ the agents must re-optimize their cost in terms of $m(\cdot,s)$ since, by   \eqref{evolucionconfeedback}, the agents move accordingly to the feedback law 
$$(x,s) \in \Omega \times [t,T] \mapsto  -\frac{\nabla u[m](x,s)}{\left|\nabla u[m](x,s)\right|}$$ rather than 
$$(x,s) \in \Omega \times [t,T] \mapsto   -\frac{\nabla u[m](x,t)}{\left|\nabla u[m](x,t)\right|}$$ (see \eqref{feedbacklaw}). In addition, their desired velocity field is re-scaled by $f(m(\cdot,\cdot))$ modelling that congestion also affects  {\it directly} the velocity of each individual agent.   Based on \eqref{evolucionconfeedback} we get that the evolution  $m$ of the initial distribution 
%
 leads, at least heuristically, to  the non-linear continuity equation 
\begin{equation}\label{continuityequationcoupled}
\begin{array}{rcl}
\partial_t m -\hbox{div}(f^2(m )\nabla u[m] \,m )&=&0, \\[6pt]
																 m(x,0) &=&m_{0}(x),
\end{array}
\end{equation} 
which is, of course, equivalent to \eqref{hughes}- \eqref{BC} (omitting the Neumann boundary condition in \eqref{BC} for $m$, which amounts to say that the trajectories followed by the agents are reflected at $\partial \Omega \setminus \mathcal{T}$).    Natural fixed point strategies  to study the existence of solutions of  \eqref{hughes} (or \eqref{continuityequationcoupled}) usually fail because of the lack of enough regularity for the solutions of both equations separately. We refer the reader to the article by \cite{di2011hughes} where an existence result is proved in the one-dimensional case $d=1$ by approximating system \eqref{hughes} by analogous systems involving  small diffusion parameters for which the well-posedness can be  shown with the help of   classical arguments in PDE theory.   Other existence results are described in \cite{amadori2012one,amadori2014existence}. 


Based on the trajectorial description presented above for both  equations in \eqref{hughes}, we consider in the next section a natural discretization of \eqref{continuityequationcoupled},  based on an Euler discretization of   equation \eqref{evolucionconfeedback} and the fact that solutions of \eqref{continuityequationcoupled} can be interpreted as the push-forward  of $m_0$   under the flow induced by \eqref{evolucionconfeedback} (cf. \cite{CS12,piccoli2011time}).

\subsection{A semi-Lagrangian scheme for a non-linear conserva- tion law}
{Equation \eqref{continuityequationcoupled} shows that \eqref{hughes} can be interpreted as a non-linear continuity equation. Note that   \eqref{value_function_depends_on_m} implies that the non-linear term $\nabla u[m]$ in \eqref{continuityequationcoupled} depends {\it non-locally} on  $m$.  In view of the previous remarks, let us describe now a SL scheme designed to numerically solve general non-linear continuity equations. The scheme we present here has been first proposed in \cite{CS12} and in \cite{CS15} in order to approximate first and second order MFGs, respectively. Then, an extension of the scheme, designed for a regularized version of \eqref{hughes}, has been implemented in \cite{carlini2016PREPRINT}. We also refer the reader to \cite{FestaTosinWolfram} where the scheme has been applied  to a non-linear continuity equation} modelling  a kinetic pedestrian model.  We recall here the scheme for the case of a two dimensional non-linear continuity equation on a bounded domain $\Omega$ with Neumann condition on  the boundary $\partial \Omega$:
 \begin{equation}\label{FP}
 \left\{
\begin{array}{ll}
 \partial_t  m + \mbox{div}( b[m](x,t)\, m ) = 0&\hspace{0.3cm} \mbox{in }\Omega \times (0,T), \\[6pt]	
  b[m](x,t) m \cdot \hat{n} =  0 \; \; &\hspace{0.3cm} \mbox{on $\partial \Omega\times (0,T)$,} \\[6pt]						 
 m(\cdot,0)= m_0(\cdot)  &\hspace{0.3cm} \mbox{in $\Omega$}.
\end{array} \right.
\end{equation}
Here, 
$b[m]: \Omega \times[0,T]\to\R$ is a given smooth vector field, depending on $m$, $m_0$ a smooth initial datum defined on $\Omega$ and
$ \hat{n}$ the unit outer normal vector to the boundary $\partial \Omega$.  Formally, at time $t\in [0,T]$ the solution of  \eqref{FP} is given, implicitly, by the image of the measure $m_0 d x $ induced by the flow $x\in \Omega \mapsto \Phi(x,0,t)$, where, given $0\leq s \leq t \leq T$, $\Phi(x,s,t)$ denotes the solution of 
\begin{equation}\label{caracteristicas_continuas}
\left\{\begin{array}{rcl}
\dot{x}(r)&=& b[m](x,r) \hspace{0.2cm} \forall \; r\in [s,T],\\[4pt]
x(s)&=& x,
\end{array}\right.
\end{equation}
at time $t$, where the trajectory is reflected  when it touches   the boundary $\partial\Omega$. 

Given $M\in \mathbb{N}$, we construct a  mesh on $\Omega$ defined by  a set of vertices $\mathcal{G}_{\Delta x}=\{x_i\in \Omega, i=1,...,M\}$ and by a set $\T$
of triangles, whose vertices  belong to  $\mathcal{G}_{\Delta x}$ and their maximum diameter is $\Delta x>0$, which form a non-overlapping coverage of $\Omega$.
We  suppose $\Omega$ to be a polygonal domain, in order to avoid  issues related to the approximation of a non-polygonal domain with triangular  meshes.\\
Given $N\in \mathbb{N}$, we define the time-step $\Delta t=T/N$ and consider a uniform  partition of $[0,T]$ given by $\{t_k=k\Delta t,\quad k=0,\hdots, N-1\}$.\\
We consider now a discretization of \eqref{FP} based on its representation formula by means of the flow $\Phi$.
For $\mu \in \R^M$, $j\in \{1,...,M\}$, $k=0,\hdots, N-1$, { we  define the discrete  characteristics as}
\begin{equation*}\label{car} 
\Phi_{j,k}[\mu]:= R(x_{j}+\Delta t \,b[\mu](x_j,t_k))
 \end{equation*}
 where $R:\mathbb{R}^2\to \Omega$ is a reflection operator, related to the Neumann boundary condition, defined as 
\begin{equation*}\label{Projection} 
 R (z):=\begin{cases}
z, &{\rm{if}}\;z\in \overline{ \Omega},\\
2\underset{w\in\Omega}{\rm{argmin}} |z-w|-z, \;\,&{\rm{if}}\;z\notin \overline{ \Omega}.
\end{cases}
\end{equation*}
We call $\{{\beta_{i}} \; ; \; i=1,...,M\}$ the set of affine shape  functions associated to the triangular mesh, such that $\beta_i(x_j)=\delta_{i,j}$  (the Kronecker symbol) and $\sum_i \beta_i(x)=1$ for each $x\in \overline\Omega$.
We define the {\em{median dual control volume}} (see \cite{Quarteroni} and \cite{Voller} for a detailed discussion on the construction of control volumes) by
\begin{align}\label{e:approxm} 
\begin{array}{c}
E_i  : = \underset{T\in \T: x_i \in \partial T}\bigcup  E_{i,T}, \; \; \; \forall i=1,\dots,M,  \\[8pt]
{\rm{where}}\quad E_{i,T}  :=\{x\in T: \beta_j(x)\leq \beta_i(x)\quad j\neq i\}. \\[6pt]
\end{array}
\end{align}


 We approximate the solution  $m$ of the problem \eqref{FP}  by a sequence $\{m_{k}\}_{k}$, where for each $k=0,\dots,N$ $m_k:\mathcal{G}_{\Delta x}\to \R$ and for each $i=1,\dots,M$, $m_{i,k}$ approximates
$$\frac{1}{|E_i|}\int_{E_i} m(x,t_k) d x,$$
where $|E_i|$ denotes the area of $E_i$. 
We compute the  sequence  $\{m_{k}\}_{k}$ by the following explicit scheme:  
\begin{equation}\label{schemefp}
\hspace{-0.4cm}\begin{array}{rl}
m_{i,k+1}&= G(m_k,i,k) \hspace{0.2cm} \forall k=0,..., N-1, \; \; i=1,...,M, \\[8pt]
m_{i,0}&=\frac{ \int_{E_{i}}m_{0}(x) d x}{|E_i|}  \hspace{0.4cm} \forall i=1,...,M,
\end{array} 
\end{equation}
where 
 $G$ is defined by
\begin{equation*}\label{definicionG}
 G (w,i,k) :=   \sum_{j=1}^M
\beta_{i} \left(\Phi_{j,k}\left[ w \right]\right)  w_j \frac{|E_j|}{|E_i|},
\end{equation*} 
for every $w\in \R^{M}$. 
\begin{remark}
In the case of a uniform standard quadrilateral mesh, the volume $E_i$ is given by $E_i=[x^1_i- \frac{1}{2} \Delta x, x^1_i+ \frac{1}{2} \Delta x] \times [x^2_i- \frac{1}{2} \Delta x, x^2_i+ \frac{1}{2} \Delta x]$, $|E_i|=(\Delta x)^2$ for each $i$ and $\{\beta_i\}$ represents the $\mathbb{Q}_1$ basis function associated to the grid.
\end{remark}
\subsection{Fast-marching semi-Lagrangian scheme for the eikonal equation}
In order to compute $\nabla u[m](x,t)$ we need to solve at each time $t$  an eikonal type equation.
For this kind  of equations,   well known and efficient techniques are Fast Marching Methods (FMM)  (\cite{Sethian}). These methods have been proposed to speed up the computation of an iterative scheme based on an upwind finite difference discretization of the eikonal equation:  the advantage is to have a complexity  $O(M \log (\sqrt{M} ))$ with respect to a complexity $O(M^2)$ of the iterative scheme. 
The FMM is a one-pass algorithm:  the main idea behind   it is that the approximation of the solution on the grid is computed following the directions given by  the characteristics equations governing the eikonal equation. This ordering allows to compute the approximated solution in only one iteration.\\
In the context of semi-Lagrangian schemes, a SL version of the FMM scheme for eikonal equations has been proposed in \cite{CristianiFalcone}, moreover a SL version of the FMM scheme  on unstructured grid has been proposed in \cite{SethianVlad2001} and \cite{CarliniFalconeHoch}.
\section{Congestion modeling}
\begin{figure}[t]
\begin{center}
\includegraphics[height=6cm]{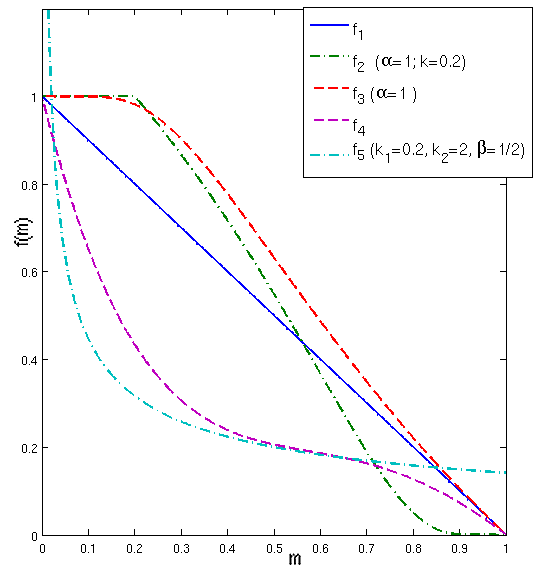} 
\caption{Fundamental diagrams. For $f_4$ the coefficients are given by $a_4=112/51$, $a_3=380/51$, $a_2=434/51$, $a_1=213/51$ and $a_0=1$.} \label{f:dia}
\end{center}
\end{figure}

The design of the function $f$ is a delicate matter that deserves a special attention. First of all we notice that the original model \eqref{hughes} is equivalent with the following system (cf. Section \ref{mod}) with the boundary conditions \eqref{BC}: 
\begin{equation*}
\left\{
\begin{array}{ll}
\partial_t m(x,t)-\hbox{div}(b[m](x,t)\,m(x,t))=0,  \\[4pt]
b[m](x,t):=f(m(x,t))\frac{\nabla u(x,t)}{|\nabla u(x,t)|}\\
|\nabla u(x,t)|=\displaystyle\frac{1}{f(m(x,t))}.
\end{array}\right.
\end{equation*}
This  system has the positive feature to avoid the numerically difficult computation of the vector field $f^2(m)\nabla u$, that, in the congested areas, is the multiplication between two quantities possibly very small $f^2(m)$ and very big (the modulus of $\nabla u$). Moreover, this formulation clarifies the role  of $f$ as quantifier of the local relation between the speed of the crowd and density. It has been proved experimentally (see e.g. \cite{seyfried2006basics, chattaraj2009comparison} and \cite{narang2015generating}) that this relation can vary as a function of several physiological and psychological factors such as the state of stress, the knowledge of the environment, etc. The graph of local density/speed is generally called \emph{fundamental diagram}, in analogy with the terminology of vehicular traffic literature. Some examples of well established choices for the diagram are:
\begin{eqnarray*}
f_1(m)&:=&1-m,\\
f_2(m)&:=&\min\left(1, \exp\left(-\alpha \frac{m-k}{1-m}\right)\right),\; \alpha>0, \; k\in(0,1), \qquad\\
f_3(m)&:=&1-\exp\left(-\alpha\frac{1-m}{m}\right), \qquad \alpha>0, \\
f_4(m)&:=& a_4m^4-a_3m^3+a_2m^2-a_1m+a_0.
\end{eqnarray*}
The choice $f_1$ appears in   \cite{hurley2015sfpe}; $f_3$ has been proposed in the early work by  \cite{weidmann1992transporttechnik} and $f_2$ can be considered as a variation of it. The diagram $f_4$ has been proposed and experimentally discussed in \cite{predtechenskii1978planning}, where the authors proposed the choice of the coefficients reported in the caption of Figure \ref{f:dia}.
In those models, the maximal speed is scaled by $1.34 m s^{-1}$, which is typically the observed maximal speed of the pedestrians;  in the same way, the maximal density before congestion (a value around $5.4 m^{-2}$) is scaled to 1.\\
Another diagram proposed by \cite{Lions}, in order to model congestion in Mean Field Games, is given by
$$f_5(m):=\frac{k_1}{(k_2 m)^\beta}, \hbox{ with }0<\beta<1/2, \hbox{ and }k_1,k_2>0. $$
In Figure \ref{f:dia}  we present some examples of the various diagrams for some adequate choices of the  coefficients.
We remark that the diagram of $f_5$ in Fig. \ref{f:dia} differs strongly from the others, since is not bounded when $m=0$ and
is not 0 when $m=1$, which means that complete congestion is not allowed.\\
The effectiveness of an approximation scheme for the system with the diagrams in Fig.  \ref{f:dia} is not obvious, since the model requires a very accurate approximation close to the  set $\{(x,t)\in \overline{\Omega} \; | \; m(x,t)=1\}$ (congestion). This brings several numerical difficulties  that should be addressed in the resolution. \\
 To the best of our knowledge, an organic comparison of different choices of fundamental diagrams in the Hughes model has never appeared in literature.

\section{Numerical simulations} 
In this section we investigate numerically the influence  of the penalty function $f$ in  \eqref{hughes}.
In order to limit instability issues in the congested zones,    in the simulations we truncate  $f$. Given $\delta>0$, in the simulations we replace $f$ by $f^\delta$ defined by
\begin{equation}\label{trunc}
f^{\delta}(m):=\max(\delta, f(m)), \hbox{ with }\ \delta\in\R_+.
\end{equation} 
Using the notations in Section \ref{discretization}, we compute the solution of the approximation of \eqref{hughes}  iteratively in the following way: given the discrete measure $m_{k}$  at time $t_k$ ({$k=0,\hdots, N-1$}), we  compute an approximation of the value function $u(\cdot,t_k)$ using the FMM scheme and the values of $m_{i,k}$ ($i=1, \hdots, M$).  This allows us to construct an approximation   $\nabla u_{i,k}$ of $\nabla u [m](x_i,t_k)$. Then, for all $i=1, \hdots, M$, we approximate $b[m](x_i,t_k)$ by $ f(m_{i,k}) \frac{\nabla u_{i,k}}{|\nabla u_{i,k}|}$ and we compute $m_{\cdot, k+1}$ using \eqref{schemefp}. 
%
%

\begin{figure}[th]
\begin{center}
\includegraphics[height=4 cm]{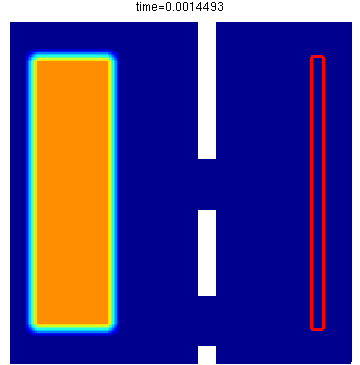} 
\includegraphics[height=4.1 cm]{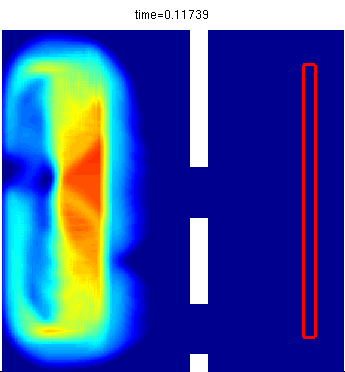} 
\includegraphics[height=3.8cm]{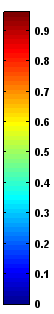} 
\includegraphics[height=4 cm]{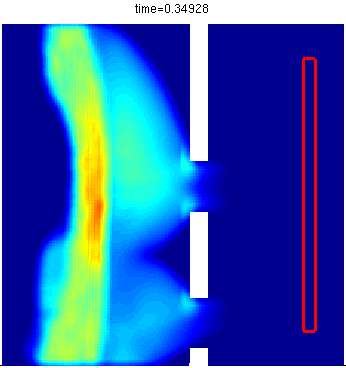} 
\includegraphics[height=4 cm]{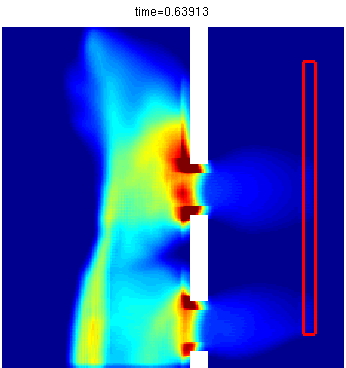} 
\includegraphics[height=3.8cm]{bar.png}
\includegraphics[height=4 cm]{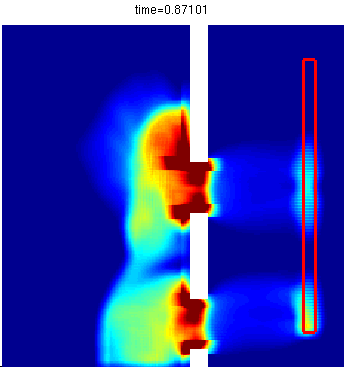} 
\includegraphics[height=4.07 cm]{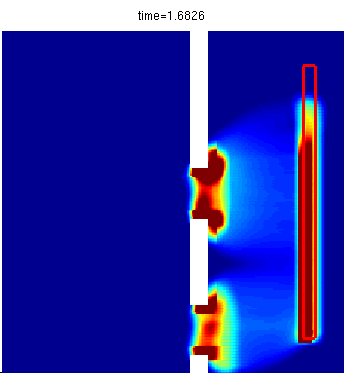} 
\includegraphics[height=3.8cm]{bar.png}
\caption{Evolution of the density with $f_1$. The red rectangle is the target region.} \label{f:t1}
\end{center}
\end{figure}

In the following simulations we fix 
$$\delta=10^{-3}, \; \; \Delta x=0.0077 \; \; \mbox{and } \; \Delta t=\Delta x/3.$$
We consider a fix scenario in a domain $\Omega:=\left([0,1]\times[0,1]\right)\setminus \Gamma$, where     $\Gamma:=\Gamma_1 \cup \Gamma_2 \cup \Gamma_3$ with
$$\begin{array}{c}  
 \Gamma_1:=[0.55,0.6]\times [0, 0.05], \; \; \Gamma_2:=[0.55,0.6]\times[0.2,0.45] ,\\[4pt]
 \Gamma_3:= [0.55,0.6]\times[0.6,1].\end{array}$$
We also fix an initial distribution given by
$$
m_0(x):=\left\{
\begin{array}{ll}
0.7 \qquad & x\in [0.1, 0.3]\times [0.1,0.9],\\[4pt]
0 & \text{ otherwise},
\end{array}\right.
$$
and a target set $\T:= [0.88, 0.92]\times [0.1, 0.95].$  
%
%

In the first test we choose $f=f_1$, which corresponds to a linear penalization of the congestion. This is  the most popular choice, see e.g. the original paper by \cite{h00} and the subsequent works by \cite{di2011hughes}  and by  \cite{carlini2016PREPRINT}. In Figure \ref{f:t1} we observe some of the basic features of the system. The mass, initially  given by   $m_0$, evolves in the direction of the target avoiding the high congested areas (top/left). For this reason the agents initially located on the extreme left hand side of the domain circumvent the center of the domain opting for less crowded regions (top/right). The density  moves towards   the  ``doors''  and the distribution of the agents takes the typical cone shape  observed experimentally in  the work by \cite{van2009pedestrian}.  Then, the crowd splits into two groups crossing the two doors. It is also possible to see (center/left-right)  as part  of the mass, initially choosing the central door, changes its strategy, preferring the bottom one,  which is  less congested. We observe that the most congested areas of the narrow passage are in contact with the walls. Note that the crowd remains congested after crossing the doors (bottom/left-right).
 This peculiar effect is due to the lack of alternative targets and, as we will see,  the choice of the model of  congestion.

\begin{figure}[th]
\begin{center}
\includegraphics[height=4.2cm]{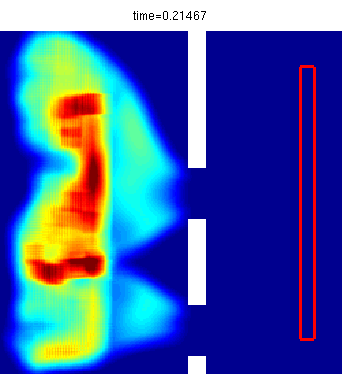} 
\includegraphics[height=4.2cm]{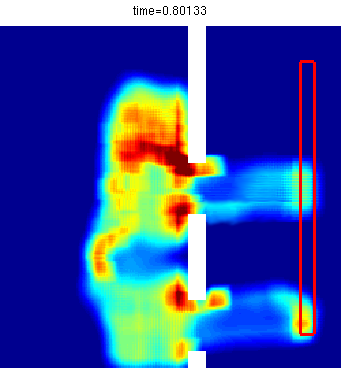} 
\includegraphics[height=3.8cm]{bar.png}
\caption{Evolution of the density with $f_2$. The red rectangle is the target region.} \label{f:t2}
\end{center}
\end{figure}

For the same initial scenario  we consider now the remaining penalty functions.
 In Figure \ref{f:t2} we display two time instants of the evolution of the system with $f=f_2$, the parameters being set as $\alpha=1$ and $k=0.2$. In this case, the behaviour of the crowd maintains some features of the previous test, however we can observe an important  difference: the maximum value reached by the density is around $0.8$ and not 1 as in the previous case. This reflects the fact that, from the perspective of the cost of each microscopic agent,   the compromise between the choices of going directly to the target and   avoiding congestion (less favorable here)  is different than in the previous case.  
%

\begin{figure}[th]
\begin{center}
\includegraphics[height=4.2cm]{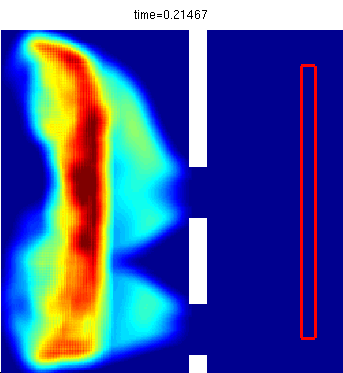} 
\includegraphics[height=4.2cm]{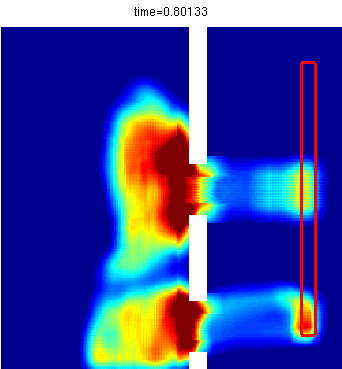} 
\includegraphics[height=3.8cm]{bar.png}
\caption{Evolution of the density with $f_3$. The red rectangle is the target region.} \label{f:t3}
\end{center}
\end{figure}

Now, we   choose $f_3$  with $\alpha=1$ as  penalty function.
In this case, we observe a macroscopic behaviour  that mixes some of the features observed so far. As before the crowd   splits into two groups associated to each door and we have that a part of the crowd change their strategies, taking into account the congestion. In this case, we observe less congestion of the crowd after crossing the doors. 
This is possibly due to a higher speed of the agents in the less crowded regions of the domain as imposed by $f_3$. 

\begin{figure}[th]
\begin{center}
\includegraphics[height=4.2cm]{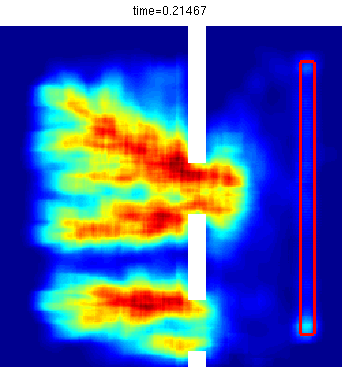} 
\includegraphics[height=4.2cm]{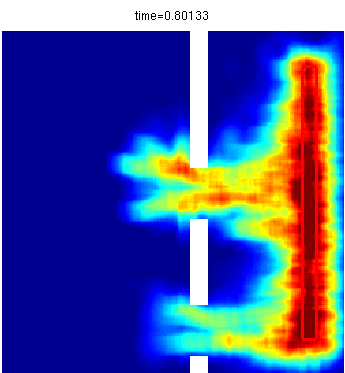} 
\includegraphics[height=3.8cm]{bar.png}
\caption{Evolution of the density with $f_4$. The red rectangle is the target region.} \label{f:t4}
\end{center}
\end{figure}

With the choice of $f=f_4$ we observe a very different behavior. In fact $f_4$ penalizes more the congested regions and it is more suitable to describe 'nervous' or 'panicked' pedestrians (\cite{predtechenskii1978planning}).  As a consequence, we do not observe   regions of high density   and the trajectories chosen by the agents result to be more  ``chaotic''.   In general, as consequence of the choice of the parameters, the time necessary to reach the target for  entire crowd is  shorter than the ones observed before. 

\begin{figure}[th]
\begin{center}
\includegraphics[height=4.2cm]{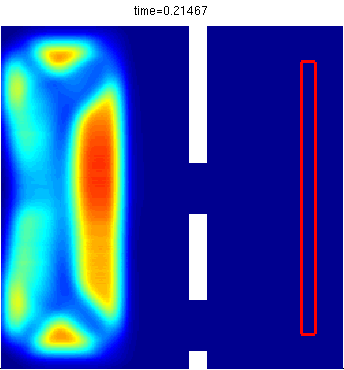} 
\includegraphics[height=4.2cm]{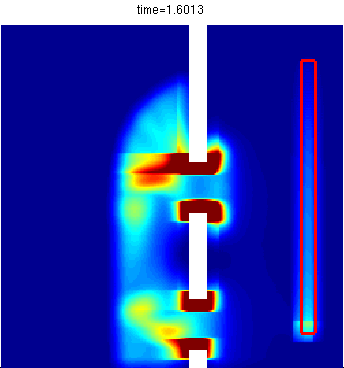} 
\includegraphics[height=3.8cm]{bar.png}
\caption{Evolution of the density with $f_5$. The red rectangle is the target region.} \label{f:t5}
\end{center}
\end{figure}

The last choice corresponds to $f=f_5$.  Here, we observe a different splitting. A small part of the crowd, moving at a quite high speed, reaches the target area in a short time. The rest of the crowd concentrates near the doors, increasing the total time to reach  the target. 

\section{Conclusions}
In this paper we have considered a   SL scheme  to solve numerically the first-order Hughes system \eqref{hughes} in the presence of various choices of the fundamental diagram $f$. The popularity of such model justifies the study of efficient and stable numerical methods for their resolution. However, many question are still open.
First of all the well-posedness of  \eqref{hughes}, and its relation with the penalty function $f$,  is not understood  in the two-dimensional  case.
 From the numerical point of view our approach requires some further work. In particular, in view of the lack of theoretical results for the continuous system \eqref{hughes}, a convergence theory for our scheme remains still as a difficult challenge.


\begin{thebibliography}{4}
\providecommand{\natexlab}[1]{#1}
\providecommand{\url}[1]{\texttt{#1}}
\providecommand{\urlprefix}{URL }
\expandafter\ifx\csname urlstyle\endcsname\relax
  \providecommand{\doi}[1]{doi:\discretionary{}{}{}#1}\else
  \providecommand{\doi}{doi:\discretionary{}{}{}\begingroup
  \urlstyle{rm}\Url}\fi

\bibitem[{Able(1956)}]{Abl:56}
Able, B. (1956).
\newblock Nucleic acid content of microscope.
\newblock \emph{Nature}, 135, 7--9.

\bibitem[{Able et~al.(1954)Able, Tagg, and Rush}]{AbTaRu:54}
Able, B., Tagg, R., and Rush, M. (1954).
\newblock Enzyme-catalyzed cellular transanimations.
\newblock In A.~Round (ed.), \emph{Advances in Enzymology}, volume~2, 125--247.
  Academic Press, New York, 3rd edition.

\bibitem[{Keohane(1958)}]{Keo:58}
Keohane, R. (1958).
\newblock \emph{Power and Interdependence: World Politics in Transitions}.
\newblock Little, Brown \& Co., Boston.

\bibitem[{Powers(1985)}]{Pow:85}
Powers, T. (1985).
\newblock Is there a way out?
\newblock \emph{Harpers}, 35--47.

\end{thebibliography}


\begin{thebibliography}{26}
	\providecommand{\natexlab}[1]{#1}
	\providecommand{\url}[1]{\texttt{#1}}
	\providecommand{\urlprefix}{URL }
	\expandafter\ifx\csname urlstyle\endcsname\relax
	\providecommand{\doi}[1]{doi:\discretionary{}{}{}#1}\else
	\providecommand{\doi}{doi:\discretionary{}{}{}\begingroup
		\urlstyle{rm}\Url}\fi
	
	\bibitem[{Amadori and Di~Francesco(2012)}]{amadori2012one}
	Amadori, D. and Di~Francesco, M. (2012).
	\newblock The one-dimensional {H}ughes model for pedestrian flow:
	{R}iemann-type solutions.
	\newblock \emph{Acta Mathematica Scientia}, 32(1), 259--280.
	
	\bibitem[{Amadori et~al.(2014)Amadori, Goatin, and
		Rosini}]{amadori2014existence}
	Amadori, D., Goatin, P., and Rosini, M.D. (2014).
	\newblock Existence results for {H}ughes' model for pedestrian flows.
	\newblock \emph{Journal of Mathematical Analysis and applications}, 420(1),
	387--406.
	
	\bibitem[{Bardi and Dolcetta(1996)}]{BardiCapuz@incollection}
	Bardi, M. and Dolcetta, I.C. (1996).
	\newblock \emph{Optimal control and viscosity solutions of
		{H}amilton-{J}acobi-{B}ellman equations}.
	\newblock Birkauser.
	
	\bibitem[{Bellomo and Dogbe(2011)}]{bellomo:2011}
	Bellomo, N. and Dogbe, C. (2011).
	\newblock On the modeling of traffic and crowds: A survey of models,
	speculations, and perspectives.
	\newblock \emph{SIAM review}, 53(3), 409--463.
	
	\bibitem[{Carlini et~al.(2013)Carlini, Falcone, and Hoch}]{CarliniFalconeHoch}
	Carlini, E., Falcone, M., and Hoch, P. (2013).
	\newblock A generalized fast marching method on unstructured triangular meshes.
	\newblock \emph{SIAM J. Numer. Anal.}, 51(6), 2999--3035.
	
	\bibitem[{Carlini et~al.(2016)Carlini, Festa, Silva, and
		Wolfram}]{carlini2016PREPRINT}
	Carlini, E., Festa, A., Silva, F.J., and Wolfram, M.T. (2016).
	\newblock A semi-{L}agrangian scheme for a modified version of the {H}ughes'
	model for pedestrian flow.
	\newblock \emph{Dynamic Games and Applications}, 1--23.
	
	\bibitem[{Carlini and Silva(2015)}]{CS15}
	Carlini, E. and Silva, F.J. (2015).
	\newblock A semi-{L}agrangian scheme for a degenerate second order mean field
	game system.
	\newblock \emph{Discrete Contin. Dyn. Syst.}, 35(9), 4269--4292.
	
	\bibitem[{Carlini and Silva(2014)}]{CS12}
	Carlini, E. and Silva, F.J. (2014).
	\newblock A fully discrete semi-{L}agrangian scheme for a first order mean
	field game problem.
	\newblock \emph{SIAM Journal on Numerical Analysis}, 52(1), 45--67.
	
	\bibitem[{Chattaraj et~al.(2009)Chattaraj, Seyfried, and
		Chakroborty}]{chattaraj2009comparison}
	Chattaraj, U., Seyfried, A., and Chakroborty, P. (2009).
	\newblock Comparison of pedestrian fundamental diagram across cultures.
	\newblock \emph{Advances in complex systems}, 12(03), 393--405.
	
	\bibitem[{Cristiani and Falcone(2007)}]{CristianiFalcone}
	Cristiani, E. and Falcone, M. (2007).
	\newblock Fast semi-{L}agrangian schemes for the eikonal equation and
	applications.
	\newblock \emph{SIAM J. Numer. Anal.}, 45(5), 1979--2011 (electronic).
	
	\bibitem[{Cristiani et~al.(2014)Cristiani, Piccoli, and
		Tosin}]{cristiani2014multiscale}
	Cristiani, E., Piccoli, B., and Tosin, A. (2014).
	\newblock \emph{Multiscale modeling of pedestrian dynamics}, volume~12.
	\newblock Springer.
	
	\bibitem[{Di~Francesco et~al.(2011)Di~Francesco, Markowich, Pietschmann, and
		Wolfram}]{di2011hughes}
	Di~Francesco, M., Markowich, P.A., Pietschmann, J.F., and Wolfram, M.T. (2011).
	\newblock On the {H}ughes' model for pedestrian flow: The one-dimensional case.
	\newblock \emph{Journal of Differential Equations}, 250(3), 1334--1362.
	
	\bibitem[{Festa et~al.(2016)Festa, Tosin, and Wolfram}]{FestaTosinWolfram}
	Festa, A., Tosin, A., and Wolfram, M.T. (2016).
	\newblock Kinetic description of collision avoidance in pedestrian crowds by
	sidestepping.
	\newblock \emph{ArXiv:1610.05056}, 1--27.
	
	\bibitem[{Hughes(2000)}]{h00}
	Hughes, R. (2000).
	\newblock The flow of large crowds of pedestrians.
	\newblock \emph{Mathematics and Computers in Simulation}, 53(4), 367--370.
	
	\bibitem[{Hurley et~al.(2015)Hurley, Gottuk, Hall~Jr, Harada, Kuligowski,
		Puchovsky, Watts~Jr, Wieczorek et~al.}]{hurley2015sfpe}
	Hurley, M.J., Gottuk, D.T., Hall~Jr, J.R., Harada, K., Kuligowski, E.D.,
	Puchovsky, M., Watts~Jr, J.M., Wieczorek, C.J., et~al. (2015).
	\newblock \emph{SFPE Handbook of fire protection engineering}.
	\newblock Springer.
	
	\bibitem[{Lions(2007-2011)}]{Lions}
	Lions, P.L. (2007-2011).
	\newblock Cours du coll\`ege de france.
	
	\bibitem[{Narang et~al.(2015)Narang, Best, Curtis, and
		Manocha}]{narang2015generating}
	Narang, S., Best, A., Curtis, S., and Manocha, D. (2015).
	\newblock Generating pedestrian trajectories consistent with the fundamental
	diagram based on physiological and psychological factors.
	\newblock \emph{PLoS one}, 10(4), e0117856.
	
	\bibitem[{Piccoli and Tosin(2011)}]{piccoli2011time}
	Piccoli, B. and Tosin, A. (2011).
	\newblock Time-evolving measures and macroscopic modeling of pedestrian flow.
	\newblock \emph{Archive for Rational Mechanics and Analysis}, 199(3), 707--738.
	
	\bibitem[{Predtechenskii and Milinski{\u\i}(1978)}]{predtechenskii1978planning}
	Predtechenskii, V. and Milinski{\u\i}, A.I. (1978).
	\newblock \emph{Planning for foot traffic flow in buildings}.
	\newblock NBS, US Department of Commerce, and the NSF, Washington, DC.
	
	\bibitem[{Quarteroni(2014)}]{Quarteroni}
	Quarteroni, A. (2014).
	\newblock \emph{Numerical models for differential problems}, volume~8 of
	\emph{MS\&A. Modeling, Simulation and Applications}.
	\newblock Springer, Milan, second edition.
	\newblock Translated from the fifth (2012) Italian edition by Silvia
	Quarteroni.
	
	\bibitem[{Sethian(1999)}]{Sethian}
	Sethian, J.A. (1999).
	\newblock Fast marching methods.
	\newblock \emph{SIAM Rev.}, 41(2), 199--235.
	
	\bibitem[{Sethian and Vladimirsky(2001)}]{SethianVlad2001}
	Sethian, J.A. and Vladimirsky, A. (2001).
	\newblock Ordered upwind methods for static {H}amilton-{J}acobi equations.
	\newblock \emph{Proc. Natl. Acad. Sci. USA}, 98(20), 11069--11074.
	
	\bibitem[{Seyfried et~al.(2006)Seyfried, Steffen, and
		Lippert}]{seyfried2006basics}
	Seyfried, A., Steffen, B., and Lippert, T. (2006).
	\newblock Basics of modelling the pedestrian flow.
	\newblock \emph{Physica A: Statistical Mechanics and its Applications}, 368(1),
	232--238.
	
	\bibitem[{Van~den Berg(2009)}]{van2009pedestrian}
	Van~den Berg, M. (2009).
	\newblock \emph{Pedestrian behaviour and its relation to doorway capacity}.
	\newblock Ph.D. thesis, TU Delft, Delft University of Technology.
	
	\bibitem[{Voller(2009)}]{Voller}
	Voller, V.R. (2009).
	\newblock \emph{Basic control volume finite element methods for fluids and
		solids}, volume~1 of \emph{IISc Research Monographs Series}.
	\newblock World Scientific Publishing Co. Pte. Ltd., Hackensack, NJ; IISc
	Press, Bangalore.
	
	\bibitem[{Weidmann(1992)}]{weidmann1992transporttechnik}
	Weidmann, U. (1992).
	\newblock \emph{Transporttechnik der fussg{\"a}nger}.
	\newblock IVT, Institut f{\"u}r Verkehrsplanung, Transporttechnik, Strassen-und
	Eisenbahnbau.
	
\end{thebibliography}
\end{document}